\documentclass[24.88pt]{article}
\usepackage{latexsym,amsfonts,euscript}

\def\pfend{\hfill{$\Box$}\vskip 0.1cm}

\title{Finite groups in which elements of the same order outside the center are conjugate
\thanks{Project supported by the NNSF of China(Grant NO. 10171074),
the NSF of Jiangsu Province College and University (Grant NO.
03KJB110112) and the Scientific Research Fund of Hunan Provincial
Education Department.}}

\author { Qian  Guohua$^1$, Shi Wujie$^2$, You Xingzhong$^{2,3}$\\
\footnotesize\small 1. Department of Mathematics, Changshu
College, Changshu, 215500, China.\\
\footnotesize\small 2. School of Mathematics, Suzhou University, Suzhou, 215006;\\
\footnotesize\small 3. Department of Mathematics, Changsha
University of Science and Technology , \\
\footnotesize\small Changsha, 410077, China. \\
}

\begin{document}
\maketitle

\begin{abstract} In this paper, we prove that if $G$ is a finite group in
which elements of the same order outside the center are conjugate,
then $G$ is abelian or is isomorphic to $S_3$, the symmetry group
of degree $3$.
\end{abstract}

\maketitle

\textbf{2000 Mathematics Subject Classification:}\,\, {\bf 20D60,
20E45}

\textbf{Keywords:}\,\, {\bf finite group, rational group, conjugacy
class
}

\section{Introduction and main results}

\quad \,\,\, In \cite{FS} and \cite{Zh}, W. Feit and G. M. Seitz,
and J. P. Zhang has solved independently the well-known  problem
posed by Syskin : If $G$ is a finite group in which elements of
the same order are conjugate, then $G\cong S_i, i=1, 2, 3$. In
this paper, we will generalize this result. Since a central
element in a group is only conjugate to itself, naturally, we
investigate finite groups in which elements of the same order
outside the center are conjugate. For convenience we call a finite
group an $OC$-group provided that all elements of the same order
outside the center are conjugate. An $OC$-group $G$ is non-trivial
if $1<Z(G)<G$. Our results are as follows.

\vskip 0.4cm {\bf Theorem A.}\,\, \em There does not exist a
non-trivial $OC$-group.\em

\vskip 0.4cm A finite group $G$ is called a rational group if
every complex character of the group is rational. Equivalently,
$g^m$ is conjugate to $g$ in $G$ whenever $g\in G$ and the integer
$m$ is relatively prime to $o(g)$, the order of $g$. Clearly, a
quotient group of a rational group is also a rational group. By
\cite{Gow} a solvable rational group is an $\{2, 3, 5\}$-group. In
order to prove Theorem A, we need investigate the rational groups
in which elements of the same odd order are conjugate and obtain
the following Theorem B. Also, Theorem B is of independent
interest, for example, we may give a concise alternative proof of
Syskin problem by Theorem B.

\vskip 0.4cm {\bf Theorem B.}\,\,  \em If $G$ is a finite rational
group in which elements of the same odd order are conjugate, then
$G/O_2(G)\cong 1, S_3, S_5, W$ or $L_3(4)<\beta>$, where $W$ is a
Frobenius group with kernel an elementary abelian group of order
$9$ and complement $Q_8$, $\beta$ is a unitary automorphism of
$L_3(4)$.\em

An interesting fact is that the quotient groups occurred in
Theorem B and the finite groups in which elements of same order
are conjugate are all the finite groups whose element orders are
consecutive integers(see \cite{BS}).

\vskip 0.4cm Groups in this paper are all finite. For a subset $A$
of a finite group $G$, denote by $\pi_e(A)$ the set of distinct
orders of elements in $A$. Other notation and terminology are
standard.
\section{Several Lemmas}

\vskip 0.4cm \quad \,\,\, {\bf Lemma 2.1.}\,\, If $G$ is a
non-trivial $OC$-group, then $|Z(G)|=2$ or $G$ is nilpotent.

\vskip 0.4cm {\bf Proof.}\,\, Suppose $G$ is non-nilpotent. Then
for any prime $p\in \pi(Z(G))$, there exists a prime $q\in
\pi(G/Z(G))$ such that $p\neq q$. Choose $x\in G-Z(G)$ such that
$o(x)=q^r$, $y\in Z(G)$ such that $o(y)=p$, Clearly, $xy,
xy^{-1}\in G-Z(G)$ and $o(xy)=o(xy^{-1}$. Hence $x^gy=xy^{-1}$ for
some $g\in G$ and so
$y^{q^r}=(x^gy)^{q^r}=(xy^{-1})^{q^r}=(y^{-1})^{q^r}$, that is,
$y^{2q^r}=1$, so $y^2=1$ and $p=2$. Thus $Z(G)$ is an elementary
abelian $2$-group. If $|Z(G)|\neq 2$, let $a_1, a_2$ be two
distinct elements in $Z(G)$, then $xa_1, xa_2$ are two elements of
the same order outside $Z(G)$ and so $x^ga_1=xa_2$ for some $g\in
G$, thus $a_1=(x^ga_1)^{q^r}=(xa_2)^{q^r}=a_2$, a contradiction.
Hence $|Z(G)|= 2$.

\vskip 0.4cm {\bf Lemma 2.2.}\,\, If $G$ is an $OC$-group, then
$G/Z(G)$ is a rational group.

\vskip 0.4cm {\bf Proof.}\,\, For any $x\in G-Z(G)$, if $m$ is an
integer prime to $o(x)$, then $x, x^m\in G-Z(G)$ and
$o(x)=o(x^m)$, so $x$ is conjugate to $x^m$ in $G$. It follows
that $\chi (x)=\chi (x^m)$ for any $\chi \in {\rm Irr}(G)$ and so
$\chi (x)$ is rational. In particular, for any $\chi \in {\rm
Irr}(G/Z(G))$, $\chi$ takes its rational value in $G-Z(G)$.
Therefore all irreducible characters of $G/Z(G)$ are rational and
so $G/Z(G)$ is a rational group.

\vskip 0.4cm {\bf Lemma 2.3.}\,\, Let $G$ be a $q$-solvable group,
$V$ a faithful irreducible $E[G]$-module over a finite field $E$
of characteristic $q$. Then there exists a faithful irreducible
character $\chi \in {\rm Irr}(G)$ such that $\chi (1)|dim_E(V)$.

\vskip 0.4cm {\bf Proof.}\,\, It is well-known that
$D=End_{E[G]}(V)$ is a finite field containing $E$, denote by
$e=|D:E|$, then $V$ is an absolutely irreducible $D[G]$-module
with degree $dim_E(V)/e$. Choose a $D$-basis of $V$ and construct
an absolutely irreducible $D$-representation $\Lambda$ of $G$. Let
$R$ be a ring of algebraic integers contained in the complex
field, $M$ the maximal ideal of $R$ such that $M\geq qR$. Then
$F=R/M$ is an algebraically closed field. Regard $D$ as a subfield
of $F$, then $\Lambda$, for emphasis, denote by $\Lambda^F$, is a
faithful irreducible $F$-representation of $G$. Let $\gamma$ be
the Brauer character corresponding to $\Lambda^F$. Because $G$ is
a $q$-solvable group, there exists $\chi \in {\rm Irr}(G)$ such
that $\chi(x)=\gamma(x)$ for any $q$-regular element $x\in G$ by
Fong-Swan theorem [\cite{F}, Theorem 10.2.1]. In particular, $\chi
(1)=\gamma(1)=dim_E(V)/e$, and so $\chi (1)|dim_E(V)$.

\vskip 0.4cm {\bf Lemma 2.4.}\,\, Let $q$ be a prime, $r, a$ two
positive integers and $r>1$. If $q^r-1|2^a \cdot 3^2 \cdot 5 \cdot
7$, then $q=2, 3, 5, 7, 11, 13, 17, 19, 29, 31, 41, 71, 127$,
where $r=2, 3, 4$ or $6$ for $q=2$, $r=2$ or $4$ for $q=3$ and
$r=2$ for other cases of $q$.

\vskip 0.4cm {\bf Proof.}\,\, If $q=2$, then $q^r-1$ is odd. Hence
$q^r-1|3^2\cdot 5\cdot 7$, so $r=2, 3, 4$ or $6$.

In the following, let $q$ be an odd prime. If $r=2^t\cdot s$,
where $s$ is an odd greater than $1$. Let
$d=(q^{2^t})^{s-1}+(q^{2^t})^{s-2}+\cdots+q^{2^t}+1$, then
$q^{{2^t}\cdot s}-1=(q^{2^t}-1)d$ and $d$ is odd. By
$q^r-1|2^a\cdot 3^2\cdot 5\cdot 7$, then $d|3^2\cdot 5\cdot 7$. A
direct checking can lead to a contradiction. Hence $r=2^t$ and by
$q^{2^t}-1=(q^{2^{t-1}}+1)(q^{2^{t-2}}+1)(q^{2^{t-3}}+1)\cdots(q^2+1)(q+1)(q-1)|3^2\cdot
5\cdot 7$, we have $t\leq 3$.

If $q=3$, by $3^r-1|2^a\cdot 3^2\cdot 5\cdot 7$ and $r=2^t$,
$t\leq 3$, then $r=2$ or $4$.

If $q>3$ and $r=2^t=4$ or $8$, then we always have
$q^4-1=(q^2+1)(q+1)(q-1)|2^a\cdot 3^2\cdot 5\cdot 7$. If
$q-1=2^c$, then $c=2$ or $c\geq 4$ and
$q^4-1=(2^{2c-1}+2^c+1)(2^{c-1}+1)2^{c+2}$, so
$(2^{2c-1}+2^c+1)(2^{c-1}+1)|3^2\cdot 5\cdot 7$, which is
impossible when $c=2$ or $c\geq 4$. Similarly, $q+1$ is not a
power of $2$. Therefore all of $q^2+1, q+1$ and $q-1$ contain odd
prime divisors which are distinct from each other, so
$p-1=2^c\cdot 3, 2^c\cdot 5, 2^c\cdot 7$ or $2^c\cdot 3^2$. If
$c=1$, then $q=7, 11$ or $19$, but in every case, $q^4-1$ does not
divide $2^a\cdot 3^2\cdot 5\cdot 7$, a contradiction; if $c=1$,
then $p+1=2\cdot 3, 2\cdot 5, 2\cdot 7$ or $2\cdot 3^2$ and so
$q=5,13$ or $17$, we can get the same contradiction.

If $q>3$ and $r=2$, then $q^2-1|2^a\cdot 3^2\cdot 5\cdot 7$. A
similar argument yields $q=5, 7, 11, 13, 17, 19, 29, 31, 41, 71,
127$ by a direct calculation.

\vskip 0.4cm {\bf Lemma 2.5.}\,\,  For any finite group $G$,
$G/Z(G)$ is not a generalized quaternion group.

\vskip 0.4cm {\bf Proof.}\,\,  see \cite{H}, p.94 Exercise 58.

\vskip0.4cm {\bf Lemma 2.6.}\,\, If $G$ is a non-trivial
$OC$-group, then $G$ is not nilpotent.

\vskip 0.4cm {\bf Proof.}\,\, Suppose $G$ is a minimal nontrivial
nilpotent $OC$-group. We infer a contradiction via several steps.

\vskip0.4cm {\bf Claim 1.}\,\, $G$ is a $2$-group.

Let $H$ be a normal $2$-complement of $G$. Clearly, $H/Z(H)$ is a
direct factor of the rational group $G/Z(G)$, so $H/Z(H)$ is a
rational group. Since a non-trivial group of odd order can not be
rational, $H=Z(H)$. Hence $G=P\times H$, where $P$ is a
non-abelian $2$-group. Clearly, $P$ is a nontrivial $OC$-group, so
$G=P$ by the minimality of $G$.

\vskip0.4cm {\bf Claim 2.}\,\, $|Z_2/Z|=2$, where $Z=Z_1=Z(G),
Z_{i+1}/Z_i=Z(G/Z_i)$.

In the following, we always assume that $exp(Z)=2^k$ and
$\overline{G}=G/Z$.

By Lemma 2.2, $G/Z$ is rational, so $Z_2/Z$ is an elementary
abelian group. Suppose $|Z_2/Z|=|\overline{Z_2}|\neq 2$. Let
$\overline{a_i}\in \overline{Z_2}$ such that $o(\overline{a_i})=2,
i=1, 2, 3$, $x\in Z$ such that $o(x)=exp(Z)=2^k$. Then $2^k\leq
o(xa_i)\leq 2^{k+1}$ and so two elements in $\{xa_1,xa_1,xa_1\}$
are of the same order, hence they are conjugate in $G$ and so
their images are conjugate in $\overline{G}$, a contradiction.
Therefore $|Z_2/Z|=2$.

\vskip0.4cm {\bf Claim 3.}\,\, All elements in $Z_2-Z$ are of
order $2^{k+1}$. Hence if $y^2\in Z$ for any $y\in G-Z_2$, then
$o(y)\leq 2^k$ and $yZ$ possesses an element of order $2^k$ in
$G$.

Suppose not all elements in $Z_2-Z$ are of order $2^{k+1}$. Since
$exp(Z)=2^k$, $Z_2-Z$ possesses an element $x$ of order $2^k$.

Let $y\in G-Z_2$ such that $y^2\in Z_2$, that is, $yZ_2$ is an
involution of $G/Z_2$. Note that $yZ_2$ possesses no elements
conjugate to $x$, then all elements in $yZ_2$ are of order
$2^{k+1}$(Otherwise, $yZ_2$ possesses an element $y_1$ such that
$o(y_1)\leq 2^k$. Let $w\in Z$ such that $o(w)=2^k$. Then
$o(y_1w)=2^k$ and $y_1w\in yZ_2$, a contradiction). Therefore $y$
is conjugate to $t$, where $t$ is some fixed element of order
$2^{k+1}$ in $Z_3-Z_2$, so all involutions of $G/Z_2$ are
conjugate. It follows that $G/Z_2$ is a cyclic or generalized
quaternion group. By Lemma 2.5, $G/Z_2$ is cyclic, then $G/Z$ is
abelian and so is $G=Z_2$, a contradiction.

\vskip0.4cm {\bf Claim 4.}\,\, $exp(Z)=2$. By claim 3, $o(x)=4$
for all $x\in Z_2-Z$ and $o(y)=2$ for all $y\in G-Z_2$ satisfying
$y^2\in Z$.

Suppose $exp(Z)=2^k\geq 4$. Let $x\in Z_2-Z$, $y\in G-Z_2$ such
that $y^2\in Z$(By Lemma 2.5, we can choose such $y$). By claim 3,
$o(xy)\leq 2^k$. So by \cite{H} Chapter 3 Lemma 1.3, we have
$$1=(xy)^{2^k}=x^{2^k}y^{2^k}[y,x]^{2^{k-1}(2^k-1)}=x^{2^k}[y,x]^{2^{k-1}(2^k-1)}
=x^{2^k}[y^2,x]^{2^{k-1}(2^k-1)}=x^{2^k}$$ a contradiction.

\vskip0.4cm {\bf Claim 5.}\,\, Let $\overline{Z_3}=Z_3/Z$. Then
$\overline{Z_3}$ is not a cyclic group.

Otherwise, by \cite{H} Chapter 3 Theorem 7.7, $\overline{G}$
possesses a cyclic normal subgroup $\overline{T}$ of index $2$.
Let $\overline{T}$ be a cyclic group of order $2^m$. Note that all
elements of orders $2^i(2\leq i\leq m+1)$ in $G$ are contained in
$T-Z$, and $G-T$ possesses no elements of order $2^{m+2}$ in $G$
(Otherwise, $G/Z$ is cyclic, a contradiction), so all elements in
$G-T$ are of order $2$. By the condition, they are conjugate in
$G$ and constitute a $G$-class with length $|G|/2$, which is
clearly impossible.

\vskip0.4cm {\bf Claim 6.}\,\, $Z_3/Z$ is an elementary abelian
group of order $4$.

First, since $G/Z_2$ is rational, $Z_3/Z_2$ is an elementary
abelian group. By claim 5, $Z_3/Z$ is not cyclic, so $G$ possesses
a subgroup $M$ such that $Z_2<M\leq Z_3$ and $M/Z$ is an
elementary abelian group of order $4$. Clearly, $M\triangleleft
G$. By claim 4, $M-Z_2$ possesses an element of order $2$. Suppose
$Z_3>M$. For all $x\in Z_3-M$, since $Z_3/Z_2$ is an elementary
abelian group, we have $x^2\in Z_2$. If $x^2\in Z$, then $o(x)=2$
by claim 4 , but $M-Z_2$ possesses the elements of order $2$, a
contradiction. It follows that $x^2\in Z_2-Z$ for any $x\in
Z_3-M$. By claim 4, all elements in $Z_3-M$ are of order $8$, so
$Z_3-M$ is a $G$-class with the length $|x^G|=|Z_3-M|\geq |M|$,
thus $|C_G(x)|\leq |G/M|$. But $xZ_2$ belongs to the center of
$G/Z_2$, we have $|C_G(x)|\geq |G/Z_2|=2|G/M|$, a contradiction.
Hence $Z_3=M$ and $Z_3/Z$ is an elementary abelian of order $4$.

\vskip0.4cm {\bf Claim 7.}\,\, The last contradiction.

Clearly, $G/Z_3$ is not cyclic. Let $xZ_3\in G/Z_3$ be any
involution. We may assume that $x$ is an element of the greatest
order in $xZ_3$, then $o(x)\geq 2$. Note that $Z_3$ contains all
elements of orders $2$ and $4$ in $G$, then $o(x)\geq 8$. Since
$exp(Z_3)=4$(see claim 6), we have $o(x)=8$. Since all such $x$
are $G$-conjugate, then all involutions in $G/Z_3$ are
$G/Z_3$-conjugate. Therefore $G/Z_3$ possesses an unique
involution and so $G/Z_3$ is a generalized quaternion group,
contradicting Lemma 2.5.

\section{Proof of theorems}

\vskip 0.4cm \quad \,\,\,\, {\bf Proof of Theorem B.}\,\, Let
$N\triangleleft G$. We prove that the conditions of Theorem A are
inherited by the quotient group $G/N$. It is clear that $G/N$ is a
rational group. Let $\overline{G}=G/N$, $\overline{x},
\overline{y}$ be two elements of the same odd order in
$\overline{G}$. We may assume that $x$ and $y$ are of the least
odd orders in $xN$ and $yN$, respectively. If $o(x)=o(y)$, then
$x$ and $y$ are $G$-conjugate and so $\overline{x}$ and
$\overline{y}$ are $\overline{G}$-conjugate. If $o(x)\neq o(y)$,
then $o(x)=p^bm$ and $o(y)=p^an$ for some odd prime $p$, where
$a<b$ and $p\not | mn$. Let $x_1=x^{mn}, y_1=y^{mn}$, then
$o(\overline{x_1})=o(\overline{y_1})$. Note that both
$x^{p^{b-a}}_1$ and $y_1$ are of odd order $p^a$, so they are
$G$-conjugate and thus $\overline{x^{p^{b-a}}_1}$ and
$\overline{y_1}$ are $\overline{G}$-conjugate. Hence
$o(\overline{x_1})=o(\overline{y_1})=o(\overline{x^{p^{b-a}}_1})$,
so $o(\overline{x_1})$ and $o(\overline{x})$ are $p'$-numbers.
This implies that $x$ is not an element of the least odd order in
$xN$, a contradiction.

By induction, we way assume that $O_2(G)=1$. Under this
hypothesis, we shall prove that $G\cong 1, S_3, S_5, W$ or
$L_3(4)<\beta>$.

\vskip0.4cm {\bf Claim 1.}\,\, $G$ possesses an unique minimal
normal subgroup.

Suppose $N_1$ and $N_2$ are two distinct minimal normal subgroups
of $G$. Let $A_1, A_2\triangleleft G$ such that
$A_1/N_1=O_2(G/N_1)$, $A_2/N_2=O_2(G/N_2)$. By induction, $G/A_1$
is isomorphic to $1, S_3, S_5, W$ or $L_3(4)<\beta>$, so is
$G/A_2$. Since $N_2$ is isomorphic to a minimal normal subgroup of
$G/N_1$ and is not an $2$-group, $N_2$ is isomorphic to a group of
order $3$ or $A_5$ or an elementary abelian group of order $3^2$
or $L_3(4)$, so is $N_1$.

If $N_2\cong A_5$, then the rational group $G/C_G(N_2)\cong S_5$.
Since elements of the same odd order in $G$ are conjugate,
$G/C_G(N_2)$ is a $\{3,5\}'$-group, and so $G/N_2$ is a
$3'$-group. By induction, $G/A_2\cong S_3, S_5, W$ or
$L_3(4)<\beta>$, then $3||G/N_2|$, a contradiction.

If $N_2\cong L_3(4)$, then the rational group $G/C_G(N_2)$ is
isomorphic to a subgroup of $Aut(L_3(4))$, by \cite{CCN},
$G/C_G(N_2)\cong L_3(4)<\beta>$. A similar argument above yields a
contradiction.

Therefore $N_2$ is of order $3$ or an elementary abelian group of
order $3^2$, so is $N_1$. But all elements of order $3$ are
contained in a $G$-class, a contradiction.

\vskip0.4cm {\bf Claim 2.}\,\, Let $N$ be the unique minimal
normal subgroup of $G$. If $N$ is unsolvable, then $G\cong S_5$ or
$L_3(4)<\beta>$.

If $N$ is unsolvable, let $N=N_1\times N_2\times \cdots \times
N_k$ is the product of $k$ isomorphic simple groups $N_i$. Let
$x\in N_1$ such $o(x)=p$, where $p$ an odd prime. Then
$x^G\subseteq N_1$. If $k>1$, let $y\in N_2$ such that $o(y)=p$,
then $o(xy)=o(x)=p$ and so $xy$ and $x$ are $G$-conjugate, thus
$xy\in N_1$, a contradiction. Therefore $k=1$ and $N$ is a
non-abelian simple.

Since $N$ is the unique minimal normal subgroup of $G$, then
$C_G(N)=1$ and so $G\leq Aut(N)$. By Theorem B in \cite{FS}, $N$
is isomorphic to one of the following simple groups:
$$ A_n, PSp_4(3), Sp_6(2), O^+_8(2), L_3(4), PSU_4(3)$$
If $N\cong A_n$, note that $A_n$ is not a rational group, then
$G\cong S_n$ or $G\leq Aut(A_6)$. If $n>6$, then $S_n$ possesses
at least two conjugacy classes of elements of order $3$, a
contradiction. If $n=6$, $G$ possesses two conjugacy classes of
elements of order $3$ by \cite{CCN}, a contradiction. Hence
$G\cong S_5$.

Again by \cite{CCN}, if $N\cong PSp_4(3), Sp_6(2), O^+_8(2)$ or
$PSU_4(3)$, $G$ always possesses two conjugacy classes of elements
of the same odd order, a contradiction; if $N\cong L_3(4)$, then
$G\cong L_3(4)<\beta>$.

\vskip0.4cm {\bf Claim 3.}\,\, In the following, we always assume
that $N$ is the unique minimal normal subgroup of odd order in $G$
and an elementary abelian group of order $q^r(q \geq 3)$, let
$A/N=O_2(G/N)$. We prove that $G$ is split in $N$, so $G=HN, H\cap
N=1$. Furthermore, $(|H|,|N|)=1$.

In fact, If $N\leq \Phi (G)$, then $q||G/N|$. Since $O_2(G)=1$,
then $O_2(G/N)=1$. So $A=N$ and by induction $G/N\cong 1, S_3,
S_5, W$ or $L_3(4)<\beta>$. Since $|N|-1$ is a $G$-class of all
elements of order $q$ in $G$, $|N|-1||G/N|$.

3.1. If $G/N\cong S_3$, then $N$ is of order $3$ and so $P\in
Sly_3(G)$ is a cyclic group of order $9$. Clearly, $P\triangleleft
G$. Therefore $P-N$ is a $G$-class of all elements of order $9$ in
$G$. Let $x\in P-N$. Then $|x^G|=|P-N|=6$ and so $|C_G(x)|=3$, a
contradiction.

3.2. If $G/N\cong S_5$, then $q=3, 5$. If $q=3$, since
$q^r-1=3^r-1|120$, then $3^r=3, 9$. Note that $G$ possesses $10$
Sylow $3$-subgroups, all $3$-elements in $G-N$ are of order $9$
and constitute a $G$-class, denote by $u^G$, then
$10(3^{r+1}-1)=|u^G|$. Therefore $|C_G(u)|=\frac{3^\cdot
4}{3^r-1}$ and so $|N|=3^r=3$, thus $|C_G(u)|=6$, contradicting
$o(u)=9$; if $q=5$, a similar argument yields a contradiction.

3.3. If $G/N\cong W$, then $q=3$ and $|N|=3, 9$. Let $H\in
Syl_2(G)$ and $P\in Syl_3(G)$. Then $H\cong Q_8$ and
$P\triangleleft G$. If $|N|=3$, since all elements of order $3$ in
$G$ are $G$-conjugate, then $G$ possesses an unique subgroup of
order $3$. Thus $P$ is cyclic, but the Sylow $3$-subgroup of $W$
is not cyclic, a contradiction. If $|N|=9$, then $|P|=81$ and
$P-N$ is a $G$-class of all elements of order $9$ in $G$. Let
$x\in P-N$. Then $|x^G|=72$ and so $|C_G(x)|=9$. In addition,
$N-1=u^G$ for some $u\in N-1$, so $|C_G(u)|=81$. It follows that
$H$ acts fixed-point-freely on $P$. Note that $2||H|$, then $P$ is
abelian and so $|C_G(x)|\geq 81$, a contradiction.

3.4. If $G/N\cong L_3(4)<\beta>$, then $q=3, 5, 7$. Let
$M\triangleleft G$ such that $M/N\cong L_3(4)$.

If $q=5$, by $q^r-1||G/N|=2^7\cdot3^2\cdot5\cdot7$, then
$|N|=q^r=5,5^2$. If $|N|=5$, then it is easy to see $N\leq Z(M)$.
Since the Schur multiplier of $L_3(4)$ is an $\{2,3\}$-group, then
$M=N \times L_3(4)$ and so all elements of order $5$ in $G$ are
not contained in a $G$-class, a contradiction. If $|N|=5^2$, since
$M/C_M(N)\cong GL(2,5)$ and $GL(2,5)$ possesses no the subgroups
that are non-abelian simple groups, then $C_M(N)=M$ and also $M=N
\times L_3(4)$, a contradiction.

If $q=7$, a similar argument above yields a contradiction.

If $q=3$, by $q^r-1||G/N|=2^7\cdot3^2\cdot5\cdot7$, then
$|N|=3^r=3, 3^2, 3^4$. If $|N|=3$, a similar argument in the case
3.3 yields a contradiction. If $|N|=3^2$, since $GL(2,3)$ is
solvable, then $N\leq Z(M)$ and so all elements of order $3$ in
$G$ are not contained in a $G$-class, a contradiction. If
$|N|=3^4$, then $|C_G(x)/N|=2^3\cdot3^2\cdot7$ for $x\in N-\{1\}$
and so $L_3(4)\cong M/N$ possesses a subgroup of order
$2^3\cdot3^2\cdot7$ or $2^2\cdot3^2\cdot7$. By \cite{CCN}, this is
impossible.

Therefore $N\cap \Phi (G)=1$ and so $G$ is split in $N$. Since
elements of same odd order in $G$ are conjugate, $N$ is a Sylow
subgroup of $G$.

\vskip0.4cm {\bf Claim 4.}\,\, $q^r-1||H|$, $H$ possesses a
faithful irreducible character $\chi $ such that $\chi (1)|r$.

Since $N-\{1\}$ is a $G$-class of all elements of order $q$ in
$G$, by claim 3, $q^r-1||H|$. Clearly, $H$ is $q$-solvable, so by
Lemma 2.3, claim 4 holds.

\vskip0.4cm {\bf Claim 5.}\,\, If $q=3$, then $G\cong S_3$ or $W$.

If $q=3$, by induction, $H/O_2(H)\cong 1, S_3, S_5, W$ or
$L_3(4)<\beta>$. Note that $H$ is a $3'$-group, so $H$ is an
$2$-group, that is, $H=O_2(H)$. By claim 4 $H$ possesses a
faithful irreducible character $\chi$ such that $\chi(1)|r$.

If $r=1$, then $\chi$ is a faithful linear character of $H$.
Therefore the rational group $H$ is of order $2$ and so $G\cong
S_3$.

If $r=2k>2$, since $q^r-1=(q^k+1)(q^k-1)||H|$, then $(q^k+1)$ and
$(q^k-1)$ are powers of $2$. But $k>1$, we get a contradiction.

If $r=2$, then the $2$-group $H$ is isomorphic to a subgroup of
$Aut(N)=GL(2,3)$. Note that the Sylow $2$-subgroup $P$ of
$GL(2,3)$ is the semidihedral group of order $16$, which is not
rational, so $H<P$. Since $8=q^r-1||H|$, then $|H|=8$ and so $G$
is a Frobenius group. Thus the rational group $H=Q_8$.

In the following we prove that $q\neq 3$ does not hold.

\vskip0.4cm {\bf Claim 6.}\,\, If $q\neq 3$, then $r=2$.

If $r=1$, then $H$ possesses a faithful linear character, so the
rational group $H$ is cyclic of order $2$, in this case, $q=3$, a
contradiction. So $q>3$ and $r>1$. By claim 4, we have $q^r-1|2^a
\cdot 3^2 \cdot 5 \cdot 7$, so by Lemma 2.4 we have $r=2$.

\vskip0.4cm {\bf Claim 7.}\,\, $q=3$ and so Theorem B holds.

Suppose $q\neq 3$. We work for an contradiction. First by claim 6,
$r=2$. Again by claim 3 we have $G=HN$ and $H\cap N=1$, $N\in
Syl_p(G)$. By induction, $H/O_2(H)\cong 1, S_3, S_5, W$ or
$L_3(4)<\beta>$. Let $\chi \in {\rm Irr }(H)$ be a faithful
irreducible character such that $\chi (1)|2$. Then $\chi (1)=2$.

If $\chi$ is not primitive, then $H$ possesses a normal abelian
subgroup $T$ of index $2$ and so $H/O_2(H)\cong S_3$. Therefore
$G$ is a solvable rational group and so $|N|=5^2$. Since $N$ is
the unique minimal subgroup of $G$, $T$ acts faithfully on $N$, in
this case, it is well-known that $C_T(y_0)=1$ for some $y_0\in
N-\{1\}$. Since $N-\{1\}$ is a $G$-class, we have $C_T(y)=1$ for
all $y\in N-\{1\}$. It follows that $T$ acts fixed-point-freely on
$N$ and so $T$ is cyclic(since $T$ is a Frobenius complement of
the Frobenius group $TN$). Note that $|y^G_0|=24$, we have
$|T|=12$ or $24$. when $|T|=24$(or $12$), the elements of order
$24$(or $12$) in $T$ are not contained in a $H$-class(since a
$H$-class of such elements is of the length $2$), so $H$ is not a
rational group, a contradiction.

If $\chi$ is primitive, by \cite{Is} p.260, $H/Z(H)\cong A_4, S_4$
or $A_5$. By induction, we have $H/Z(H)\cong S_4$. Note that $H$
is rational and $Z(H)$ is cyclic, then $|Z(H)|\leq 2$. Since
$\chi$ is primitive, $Z(H)>1$, so $H\cong GL(2,3)$. But by
\cite{Hu} E11.8, $GL(2,3)$ is not a rational group, a
contradiction.

\vskip 0.4cm As a corollary of Theorem B, we give a concise
alternative proof of Syskin problem.

\vskip 0.4cm {\bf Corollary 3.1.}\,\, \em If $G$ is a finite group
in which all elements of the same order are conjugate, then
$G\cong S_i, i=1,2,3.$\em

\vskip 0.4cm {\bf Proof.}\,\, Note that the condition of corollary
3.1 is inherited by the quotient groups and $G$ is a rational
group, by Theorem A, $G/O_2(G)\cong 1, S_3, S_5, W$ or
$L_3(4)<\beta>$. Note that $S_5$ possesses two $G$-classes of
elements of order $2$ and $W$ and $L_3(4)<\beta>$ possess three
$G$-classes of elements of order $4$, so $G/O_2(G)\cong 1$ or
$S_3$. If $G=O_2(G)$, since all elements of the same order are
conjugate, if $G\neq 1$, then $G$ possesses an unique element of
order $2$. Therefore $G$ is a cyclic or generalized quaternion
group and so $G\cong S_2$. If $G/O_2(G)\cong S_3$, let
$exp(O_2(G))=2^k$. Then $O_2(G)$ contains all elements of orders
$2^i(1\leq i \leq k)$ in $G$ and so $o(x)=2^{k+1}$ for all
$2$-elements $x\in G-O_2(G)$. It follows that all $2$-elements in
$G-O_2(G)$ constitute a $G$-class with length
$|x^G|=3|O_2(G)|=|G|/2$ and so $|C_G(x)|=2$. Hence $k=0,
O_2(G)=1$, that is, $G\cong S_3$.

\vskip 0.4cm {\bf Proof of Theorem A.}\,\, Suppose that $G$ is a
non-trivial $OC$-group. We work for a contradiction.

By Lemma 2.1 and Lemma 2.6, $G$ is non-nilpotent and $|Z(G)|=2$,
so $G$ is a rational group satisfying the conditions of Theorem B,
we have $G/O_2(G)\cong S_3, S_5, W$ or $L_3(4)<\beta>$. Let
$M\triangleleft G$ such that $O_2(G)<M<G$ and $|G/M|=2$.

Let $x\in G-M$ be not an $2$-element, then $x=yz, yz=zy$, where
$y$ is an $2$-element and $z$ is of odd order. It is easy to see
that $y\in G-M$ and $o(z)$ is a fixed positive integer, for
example, $o(z)$=3 if $G/O_2(G)\cong S_5$ or $L_3(4)<\beta>$; $z$
does not exist if $G/O_2(G)\cong S_3$ or $W$, that is, the
elements in $G-M$ are all $2$-elements.

Let $P\in Syl_2(G)$, $P_1=M\cap P\in Syl_2(M)$, $exp(P_1)=2^k$. By
the condition, there do not exist $2$-elements $x\in G-M$ and
$y\in M-Z(G)$ such that $o(x)=o(y)$. Note that $G-M$ is not a
$G$-class, otherwise $G$ is a Frobenius group with kernel $M$,
contradicting $Z(G)>1$. We consider the following two cases.

Case 1, $P-P_1$ possesses no involutions.

For any $x\in G-M$, if $x$ is an $2$-element, then $o(x)=2^{k+1}$;
if $x$ is not an $2$-element, then $o(x)=3\cdot 2^{k+1}$. Hence
$G-M$ possesses just two $G$-classes of elements of orders
$2^{k+1}$ and $3\cdot 2^{k+1}$, respectively. It follows that
$|C_G(x)|=4$ for any $x\in G-M$, contradicting that $G-M$
possesses the elements which are not $2$-elements.

Case 2, $P-P_1$ possesses an involution.

In this case, $P_1-Z(G)$ possesses no involutions, so $P_1$ is a
cyclic or generalized quaternion group.

If $P_1$ is cyclic, then $G/O_2(G)$ is neither isomorphic to $S_5$
nor $L_3(4)<\beta>$. When $G/O_2(G)\cong S_3$ or $W$,
$\pi_e(G-M)=\{2, 2^{k+1}\}$ and $G-M$ possesses just two
$G$-classes. Hence $P$ is cyclic, a contradiction.

If $P_1$ is a generalized quaternion group.

(a). If $G/O_2(G)\cong S_3$ or $W$, then $G-M$ possesses just two
$G$-classes and $\pi_e(G-M)=\{2, 2^{k+1}\}$. It follows that
$|C_G(x)|=4$ for any $x\in G-M$ and so $2^{k+1}=4, k=1$. Therefore
$P_1$ is an elementary abelian group, a contradiction.

(b). If $G/O_2(G)\cong S_5$, since $P_1$ is a generalized
quaternion group, then $O_2(G)$ is a cyclic subgroup of $P_1$ of
index $4$ and $C_M(O_2(G))> O_2(G)$, and so $O_2(G)=Z(M)$.
Clearly, $M=M'$(Otherwise $P_1$ could not be a generalized
quaternion group), hence $O_2(G)$ is the Schur multiplier of
$A_5$, where $A_5\cong M/O_2(G)$, so $O_2(G)=Z(G)$. Therefore we
have $M=SL_2(5)$ and $G=2.A_5.2$, where $2.A_5.2$ is the Schur
covering group of $A_5$, but in this case, both $G-M$ and $M$
possess the elements of the order $6$, a contradiction.

(c). If $G/O_2(G)\cong L_3(4)<\beta>$, then $M/O_2(G)\cong
L_3(4)$. Since $P_1$ is a generalized quaternion group, we can
easily get a contradiction by investigating the Sylow
$2$-subgroups of $L_3(4)$.

For solving Syskin's problem we all used the classification of
finite simple groups directly or indirectly in \cite{FS},
\cite{Zh} and Theorem B of this paper.

{\bf Problem} Whether or not exist a proof for solving Syskin's
problem which does not depend on the classification of finite
simple groups?

\vskip 0.4cm{\bf Acknowledge} The authors would very much like to
thank Professor J. P. Zhang for his helpful comments at the stage
of preparation of this paper. In fact, he pointed out the result
of Lemma 2.1. \pfend

\end{document}